\documentclass[12pt]{article}
\usepackage{amssymb}
\def\be{\begin{equation}}
\def\ee{\end{equation}}

\def\R{{\sf I\kern-.2em R}}
\def\N{{\sf I\kern-.2em N}}
\def\C{\kern.1em{\raise.47ex\hbox{$\scriptscriptstyle
$}}\kern-.40em{\sf C}}
\def\Z{{\sf Z\kern-.32em Z}}

\def\hat{\widehat}

\def\hat{\widehat}
\def\be{\begin{equation}}
\def\ee{\end{equation}}

\newtheorem{theorem}{\noindent Theorem}
\newtheorem{lemma}{\noindent Lemma}
\newtheorem{definition}{\noindent Definition}
\newtheorem{corollary}{\noindent Corollary}

\newtheorem{remark}{\noindent Remark}

\begin{document}

\begin{titlepage}

\begin{center}
{\LARGE \bf GEOMETRY AND DYNAMICS ON THE FREE SOLVABLE GROUPS}

\bigskip
\bigskip

{\large A.~M.~VERSHIK $^{1}$}
\bigskip

{\sl
$^{1}$ Steklov Institute of Mathematics at St.~Petersburg,\\  
Fontanka 27,
119011 St.~Petersburg, Russia.\\
Partially supported by RFBR, grant 99-01-00098.} 
\end{center}

\centerline{15.06.2000.}

\bigskip

\vspace{0.3in}

\begin{abstract}

In this paper we give a geometric realization of free solvable groups, 
and study its Poisson-Furstenberg boundaries, we also discuss
the construction of normal forms in the solvable groups.

\end{abstract}

\vspace{0.3in}
Keywords:free solvable groups, homology, fundamental 2-cycle,
boundary, normal form.

\end{titlepage}

\newpage

\section{Introduction}

Free solvable groups were studied by algebraists in the 
40s--50s in works by F.Hall, W.Magnus and others \cite{Hall,Mag,Schm}).
The main results were concerned with imbedding to the wreath product.
Later in 60-s the growth of the lower central series \cite{Sok}, and so
called Golod--Shafarevich series were discussed in the literature 
(e.g.~see \cite{Petr}) were considered. Now these groups have
become an object of study from the viewpoint of harmonic analysis and
asymptotic characteristics. For this we need in more precise model of these
groups. Unfortunately the imbedding to the wreath product which was stuied
before is not effective because the image of the groups is difficult 
 to describe explicitly.

The free solvable groups of level two and higher have not exact matrix
representations,in a sense they are "infinite-dimensional" (or "big") groups,
as we shall see further  -- in contrary, say,  to the free nilpotent groups. In
this paper we give a {\it new topological model of the free metabelian
groups}, i.e.~the free solvable groups of level two, and calculate their
boundaries. Our realization  differs from the known ones (Nilsen--Schraer
basis of commutant) by its invariance.

Some of wreath products are very similar to the free solvable groups, its
had been considered for a long time in the theory of growth and random walks
as a natural source of examples and counterexamples. For example in
\cite{VSel} the wreath product $\Z \wr \Z$ was used as an example of a group
with superexponential growth of Folner sets. Later in \cite{KV} the group
$\Z^d \wr \Z$, $d > 2$ --- were used as examples of a solvable and thus
amenable groups with non-trivial boundary. The term ``wreath product'' was not
used in these works,  but the suggested terms ``lamplighter group'' and
``group of dynamical configurations'' became popular. Thus the wreath products
$\Z^d \wr \Z$ for $d > 2$ give an important example of amenable groups with
positive  entropy - it was a surprise --- a common opinion before was
that all amenable groups have zero entropy (\cite{Av}). The
Furstenberg-Poisson boundary in this example was not calculated completely,
only a natural candidate was presented --- moreover,  it is this candidate
that was used to prove positivity of the entropy  (owing to the entropy
criterion \cite{KV}). These kind of  examples were used also in the theory 
 of index of von Neumann factors (see \cite{Hag}).  Recently, the wreath
products were used by a student A.~Dyubina - see \cite{Dyub} in  construction
of an example of quasi-isometric groups, one of which is solvable, and the 
other one is not virtually solvable, as well as  in construction of an example
of intermediate drift growth. But the boundary was  not explicitly calculated
upto now - we give the precise description as a corollary of the theorem about
the boundary for free solvable group.

Now it becomes clear that the {\it free solvable groups are much more natural
and important class of examples than the wreath products with lattices, and the
effects which were discovered in the wreath products manifest themselves yet
more explicitly in the free solvable groups.}

In this paper we shall give: 

\smallskip

a new topological model of metaabelean
(= free solvable group of level 2), 

\smallskip

a new (geometrical) normal form for the elements of it

\smallskip

and finally  we describe 
the Poisson-Furstenberg boundary of these groups. 
\bigskip

Using the same method and reduction to the free solvable group, we  also descirbe  
the boundary of wreath products - old problems which
appeared in  \cite{KV}.
 \newpage
\section{Topological model of free solvable groups of the level two}

We start with a new, as far as we know, model of the free solvable group of level two
having a topological interpretation.

Let $Sol^2_d$ be the free solvable group of level two with $d$ generators,
i.e.~the universal object of the variety of solvable groups of level two with $d$
generators. This group may be defined in a more constructive way as the factor group
of the free group with $d$ generators with respect to the second commutant:

$$ Sol^2_d=F_d/(F_d)''$$
where $F_d$ is free groups with $d$ generators, and $G''$ is the second
comutant of the group $G$. Sometimes we omit "2" in $Sol^2_d$ because 
in this paper we consider solvable group of level two only.

The commutant of this group
is an abelian group with infinite number of generators which may be
indexed by the elements of the $d$-dimensional lattice. Namely, they are the images
of the commutators $[x_i,x_j]$, $i,j=1, \dots ,r$;
$i \ne j$ of the original generators under the action of the inner automorphisms
defined by the monoms $x_1^{k_1} \dot x_2^{k_2} \dots x_r^{k_r}$ (the order of factors
does not matter since such automorphism of the commutant, as one can easily check,
depends only on the degrees of variables $k_1,k_2, \dots,
k_d$ and does not depend on their order). The group $Sol_d$ is isomorphic
to the extension of the commutant by the group ${\Z}^d$ acting on the commutant and
some $2$-cocycle. (The author follows the terminology in which the extending group is the
group which acts by automorphisms rather than the group 
on which it acts, as is accepted in algebra).

For $d=2$, the commutant is completely defined, since the described generators 
of the commutant are not subject to any relations, and {\it they are  the
free generators of an abelian group}. But for $d \geq 3$ the description
becomes more cumbersome, since the generators are subject to some relations.
One can easily see this:  for example the element
$$
x_1x_2x_3x_1^{-1}x_2^{-1}x_3^{-1}
$$ 
admits two different notations in pairwise commutators of the
group generators and adjoined elements,
i.e.~in the above described generators of the commutant. Easy to see that this element
is a cycle on the one-dimensional skeleton of the three-dimensional cube, and
there are two ways to
decompose it into the product of cycles of plane faces.
Thus the use of these generators is not convenient. It was proved in classical works
(\cite{Hall} and others) that this group is embedded into the wreath product
$\Z^d \wr \Z^d$, but this embedding is not fit for our purposes either, 
since the image of these 
groups under this embedding is not easy to describe.

However, there is a direct way to describe the commutant and the whole group which
is partially borrowed from analogies with the {\it theory of cubical
homologies, see f.e. \cite{Hil}}.  
The same method also allows to describe easily the cocycle defining the extension.

Let us consider the lattice as a one-dimensional complex, i.e.~a topological space, namely,
as the union of all shifts of the coordinate axes by integer vectors. In
another word this is the Caley graph of $\Z^d$ under the ordinary generators as
one-dimensional complex. To distinguish it from the lattice as a discrete
group, denote this one-dimensional complex by $E^d$. Consider the additive
group ${\bf B}_d$ of oriented closed $1$-cycles on the space $E^d$ as a
one-dimensional complex. We will consider nontrivial in homotopy sense
cycles or simply the groups of the first homologies of $E^d$ with integer
coefficients - $H_1(E^d)\equiv B_d/Z_d$, where $Z_d$ is the group of the
homotopivcally trivial cycles. The generators  of $H_1(E^d)$ are {\it
elementary cycles (plackets in the physical terminology, or standard
$1$-cycles  in sense of the theory of cubical homologies, see\cite{Hil})},
i.e.~the cycles that go around two-dimensional cell of the lattice.  A path
around the cell with nodes $(0,e_i,e_i+e_j,e_j)$ in the indicated order, where
$e_i$, $e_j$ are coordinate unit vectors, is called the {\it elementary
$(i,j)$-coordinate cycle}. There is a natural action of the group ${\Z}^d$ by
shifts on the group of homologies $H_1(E^d)$. Each elementary cycle is a
translation of one of the elementary coordinate cycles.

 {\it Now we define a $2$-cocycle $\beta(\cdot,\cdot)$ of the group
 ${\Z}^d$ with values in the group of homology $H_1(E^d)$ or in the group
of cycles ${\bf B}_d$}, and then the element of the group $H^2(\Z^d; H_1(E^d))$
with respect to action defined above. More exactly, we define  at once the
cohomology class of cocycles. For this, we associate with every element $v \in
{\Z}^d$ of the lattice an arbitrary connected path $\tau_v$ connecting the
zero with the element $v$ (``path'' as well as ``lattice'' are understood
literally: a path is a continuous mapping of the half-line $R_{+}$ or a
segment into the lattice as a topological space which sends integer points
$\N\in R_{+}$ to integer vectors, and which is linear at each integer
segment). Then, given a pair $(v,w)$ of elements of the lattice $\Z^d$, we
define a {\it cycle} $\beta(v,w) \in {\bf B}_d$ formed by three paths: $$ 
(\tau_v,v+\tau_w, -\tau_{v+w}). $$  
This cycle regarded as an element of
 ${\bf B}_d$ is exactly the value of the $2$-cocycle $\beta (v,w)$,
or more exactly as element of the group $H_1(E^d)$.

\begin{lemma}
Different choices of the system of paths  $v \mapsto  \tau_v$ lead to cohomological
cocycles.
\end{lemma}
\begin{proof}
{
Indeed, if $v \mapsto \tau_v$ and $v \mapsto \rho_v$ are two such systems, then the cycle
$(v \mapsto \tau_v ,(v \mapsto\rho_v)^{-1})$ 
realizes the cohomology.
}
\end{proof}

As a corollary we obtain that the cocycle correctly defined the class of
cohomology,  Denote by $\bar \beta$ the cohomology class of the constructed
cocycle. Note that by construction this cocycle is trivial in the group of
paths, however it is non-trivial as a cocycle in the group of cycles (or
homologies) of $E^d$. So we obtain the extension of infinitely generated
abelian group $H_1(E^d)$ by the group $\Z^d$ with the 2-cocycle $\beta$. The
generators of this group are the genrator of $\Z^d$.
  \begin{theorem} 
The extension of the group $H_1(E^d)$ by the group of shifts
${\Z}^d$ with the cohomology class $\bar \beta$ is canonically isomorphic
to the free solvable group of level two with $d$ generators $Sol_d$. 
The image of the group  $H_1(E^d)$  under this isomorphism is exactly
the commutant of the group $Sol_d$, the image of the elementary
$(i,j)$-coordinate cycle being the commutator of the generators $[x_i,x_j]$,
and the action of $\Z^d$  on cycles turning 
into the action of $\Z^d$ by the inner automorphisms on the commutant of the
group $Sol_d$. 
\end{theorem}
\begin{proof}
{Let us define homomorphism from $Sol_d$ to the groups we have defined by
putting the generators of $Sol_d$ to the generator of $\Z^d$. It is clear
that the commutators if the elements of $\Z^d$ go to the subgorup of cycles
(homologies), so we have the surjection of $Sol_d$ onto our group. It is
clear also that the kernal is trivial.}
\end{proof}

\smallskip

\begin{remark}
{It is possible to prove that
$H^2({\Z}^r; {\bf B}_r)=\Z$,  and the constructed cocycle is a generator of
this group.}
\end{remark}

This construction admits far generalizations. For example, it becomes clear how to
understand a continual analogue of the free solvable groups of level two --- one should
replace the group of cycles on the lattice by an additive group of some
$1$-cycles on $R^d$ (or De Rham flows), and define the action of  $R^d$ and the
cocycle exactly as above. Such generalization makes clear the above remark that
 one should consider the free solvable group as
an infinite-dimensional group.  

More exactly.  Let $M$ a homogeneous space of the Lie group $G$ and
$H_1(M)$ is the first homologies with compact support with scalar coefficients
(say $\C$). The group $G$ acts on $H_1(M)$, and we can define the cocycle
$\beta$ in the same way: to choose the path $\gamma_x$ from some fixed
point $0$ to an arbitrary point $x$ which depends continiously on the
point $x$ in the space of pathes, and then define the 2-cocycle with the same
formula. 

It is interesting to define also the analogue of this
construction with smooth differenctial 1-forms instead of $H_1$.
andto define canonical 2-cocycle and element of $H^2(G;\Omega^1(M)$
related to de Rham cohomology.  The same technique
may be applied to construct the free solvable groups of higher levels --- for example,
the group of level three may be represented in the same way, since its commutant
is the free solvable group with infinite number of generators, hence 
it is natural to represent the second
commutant as the group of $2$-cycles on the lattice, and by induction we
can constract 
$$\Z^d \ltimes_{\beta}(H_1(E^d) \ltimes_{\theta}H_1 (H_1 (E^d)))$$

where $H_1(H_1(E^d))$ is the group of the first homologies on the abelian
group $H_1(E^d)$ with 
respect to ${\it free}$ generators, and $\theta$ the 
corresponding 2-cocycle  of the same type as $\beta$ above.

\newpage
\section{Space of the pathes and normal forms}

We are going to present a general simple technique which reduces the word 
identity problem in
finitely generated groups to the combinatorial geometry on the lattice. We shall see
that the word identity problem has a pronounced geometric character.
The word identity in a group is equivalent to a notion of equivalence of paths 
on the lattice which depends on the group; to find a normal form is to solve an
isoperimetric problem, etc. In fact, we replace the space of paths on the Cayley graph
by a canonically isomorphic space of paths on the lattice. In some cases like
the one we study below this method is very efficient.

\smallskip
Given a group $G$ with a system $S=S^{-1}$ of $d$ generators, each its element
can be represented by a word in the alphabet $S$, and with this word we may
associate a connected path on the lattice ${\bf \Z^d}$ starting at zero as
follows. Identify the $i$th generator with the $i$th coordinate unit vector of
a fixed basis of the lattice, $i=1, \dots, d$. Now associate with each {\it
oriented\/} edge of the lattice a generator or its inverse which corresponds
to the coordinate axis parallel to this edge, taking the generator, if we pass
the edge in the direction of growing distance from zero, and the inverse
generator, if the direction is opposite. It is clear, that the space of pathes
of the given l (finite or infinite length) in the Cayley graph  of the groups
with the fixed set of $d$ generators  canonically isomorphic to the space of
pathes of the same length on the lattices $\Z^d$.

 Thus each word is a path on the lattice: the empty
word (the unity of the group) turns into the path consisting of one zero point
of the lattice, let us call it trivial. Adding some generator (or its inverse)
to the end of a given word means adding the oriented edge corresponding to this
generator or its inverse, depending on orientation,  to the end of the
constructed path.   

\begin{definition}
{
Two paths are $G$-equivalent if they define the same element of the group $G$.
}
\end{definition}

Thus the word identity problem is reformulated as the problem of $G$-equivalence
of paths on the lattice. Of course this reformulation looks like tauthology,
but sometimes it is very useful.

If the relations between generators of the group are generated by the elements of
the commutant of the free group, then it suffices to define the equivalence of a closed
path (=cycle) to the trivial path (=cycle). This is the case in some of the below examples.
But sometimes (e.g.~the Heisenberg group) a non-closed path may be equivalent to a closed
one. In this case it does not suffice to define the equivalence of closed paths to the
trivial one.

\bigskip

{\Large Example 1.}

{\it Free abelian groups}
\smallskip
 If $G= \bf \Z^d$, 
then two paths are $G$-equivalent if their ends coincide. The group is identified
with the lattice. Every closed path (terminating at zero) is equivalent to the
trivial one and defines the unity of the group.

\bigskip
{\Large Example 2.}

{\it Free groups}

\smallskip
Another trivial example is given by the free group. In this case two paths are equivalent 
if they will coincide if we 
successively cancel in each path all {\it neighbouring} edges
differing only by orientation.

\bigskip
{\Large Example 3.}

\smallskip
{\it Free nilpotent groups of level\/} 2.

\smallskip
A less trivial example. Let $G$ be the free nilpotent group of level 
$2$ with $d$ generators. First consider the case $d=2$ --- this is the discrete
Heisenberg group, i.e.~the group of integer-valued upper triangular
matrices of third order with
ones on the main diagonal. Denote the generators by
$a$, $b$, then the relations are as follows:
$$
[a,b]a=a[a,b],\quad [a,b]b=b[a,b],\quad [a,b]=aba^{-1}b^{-1}.
$$

In this case a closed path is equivalent to the zero path if and only if the 
algebraic (oriented) area enclosed by it is zero. And two paths
$\gamma_1$  and  $\gamma_2$  are equivalent if the closed path formed by
$\gamma_1 \gamma_2^{-1}$ is equivalent to the trivial path.

The same conclusion is also true for the continuous Heisenberg group --- the area is
exactly the value of the symplectic $2$-form (defining the group) on the corresponding
$2$-cycle. This fact is known for a long time 
and is used in symplectic geometry and the 
control theory.

For $d>2$, the equivalence of a closed path to the zero one is described by the same
criterion but applied to all projections of the closed path on the two-dimensional
coordinate subspaces of the lattice.

The following key example seems to be new.

\bigskip
{\Large Example 4.}

\smallskip
{\it Free solvable groups of level two - free metaabelian groups}.

Let $G=Sol_d$. By definition, the relations between generators lie in the
second (and hence in the first) commutant of the free group. Thus it suffices
to define the equivalence of closed cycles to the trivial cycle. But we give
the  equivalence condition  for arbitrary paths.

Given a path $\gamma$, an edge $\rho$ of the lattice may occur in $\gamma$
with different orientation $+,-$. Denote by $\gamma(\rho)$ the algebraic sum
(including signs) of all occurrences of this edge in the path $\gamma$.

\begin{lemma}
{
 Two finite paths $\gamma_1$, $\gamma_2$ on the lattice $\bf \Z^d$ 
are $Sol_d$-equivalent if and only if $\gamma_1(\rho)=\gamma_2(\rho)$
for every edge $\rho$.
}
\end{lemma}

\bigskip

In other words, the equivalence means that all edges of the lattice occur in both paths
with equal multiplicities (taking into
account the directions) independently on the order. Recall that
in Example~2 (the free group) only neighbouring edges of opposite orientation were
canceled.

\begin{proof}
{The proof follows immediately from the description of the commutant of the group
$Sol_d$. Indeed, since all commutators of $Sol_d$ commute, this means that
two paths differing by the order of plackets occurring in them are equivalent.
Thus the equivalence class depends only on the total multiplicity of 
edge taking into account orientations.
}
\end{proof}

This criterion may be easily reformulated as a solution of the word identity problem
in the group $Sol_d$ in inner
terms of the words, i.e.~as a normal form of group elements,
but the above criterion is more useful for the sequel.
The detale description of this normal form was investigated by student
S.Dobrunin. The main preference of this group is the cancellation property for edges
- in order to make a label on the given edge we can restrict ourself with the
visits of this edge only (and ignore other parts of path). This is not the
case  for the solvable group of the levels more than two.
A nice question -to find the groups with the natural geometrical equvalence of
the pathes. The language of the equivalence of the words in initial alphbet
could be more cumbersome than the language of the equvalence of the pathes (but
of course both are coinsided), The group $Sol_d$ is just an example of that
effect.
 Our interpretation of the classes of the group could be considered as 
"normal form" of the word - the difference with ordinary normal form is the
following: we are givong the {\it invarianats of the classes - function on the
edges - instead of giving some concrete representor of classes.} Of course it
is not difficult
using our method to give the normal form in usualk sense, but this form won't
be minimal.

\newpage 
\section{The boundary of the free solvable groups and of the wreath products}

Now consider the space of infinite paths on the lattice ${\Z}^d$ beginning at zero.
It may be identified with the space $S^{\N}$ of infinite sequences
in the alphabet $S$ of generators and their inverses. We provide it with 
the natural topology and the Bernoulli measure $\mu^{\infty}$ with equal
probability on the generatos and its inverses.

Consider also the space $F(E^d;\Z\cup +\infty)\equiv F^d$ of all functions with
integer or infinite values on the set of edges of the lattice $E^d$, and
introduce the mapping $\Phi$ of the space of infinite paths
$S^{\N}$ in $F^d$ which associates with an infinite path
$\gamma$ and an edge $\rho$ the number $\lim_{n \to \infty}\gamma_n(\rho)$,
if the limit exists (i.e.~stabilizes), and $+\infty$ otherwise. Here 
$\gamma_n$ is the initial segment of the path $\gamma$.

We will consider a "simple  symmetric" random walk on the free solvable group 
$Sol^d$ of the level 2 with $d$ generators wich means that the initial
measure (transition probability) is a uniform measure with charge $(2d)^{-1}$
on the canonical generators and its inverses. We will describe 
{\it the Furstenberg- Poisson boundary - $\Gamma(Sol^d, \hat\mu)$} (see for
drfininiton \cite{KV} and \cite{V}) of that random walk, using a reduction
to the classical simple symmetric random walk on the group $\Z^d$. Namely, each
path (trajectory) in the $Sol^d$ can be projected to the $\Z^d$ so, we have
a map from space of pathes in $Sol^d$ to the space of pathes  $\Z^d$. We
will use the map very intensively.

Recall a fundamental fact from the theory of simple random walks on the
lattices: if $d=1,2$, then the random walk is recurrent, i.e.~with probability
one it passes infinitely many times through any node and any edge, and if
$d>2$, the random walk  is not recurrent, in particular, each edge with
probability one occurs in a trajectory of the walk with finite multiplicity
(see \cite{Spi}). Further we shall consider this case $d>2$. Let us call   the
functions from $F^d$ stable functions.  (We exclude infinite values, since 
the walk is non-recurrent, thus the values of stable functions are finite.)
Two infinite paths with the same image under $\Phi$ (finite or not) are called
stable equivalent. Thus stable equivalence coincides with the equivalence
considered in the previous section. Our aim is to prove that this equivalence
coincides with the equivalence defined by the boundary partition $\sigma$ on
$S^{\N}$, see section~3. In other words, we have to prove that there is no
measurable function (it suffices to consider functions from $L^2$) that is
orthogonal to all stable functions and is $\sigma$-measurable, i.e.~is
constant on the paths whose finite segments are equivalent in the above sense.

\begin{theorem}
{
1. For $d \leq 2$, the space of classes of infinite stable equivalent paths,
and the boundary $(\Gamma(Sol_d, \mu), \mu_\sigma)$, where $\mu$ is the uniform measure
on generators of the group $Sol_d$, are trivial ($\bmod0$), 
i.e.~consist of a single point. More
exactly, the mapping $\Phi$ sends almost all paths to the identically infinite function
on edges.

2. For $d>2$, the space of classes of stable equivalent paths is not trivial,
and it is canonically isomorphic to the boundary $(\Gamma(Sol_d,
\mu), \hat \mu)$. More exactly, the mapping 
$$
\Phi: (S^{\N}, \mu^{\infty}) \to (F^d, \hat \mu)
$$ 
is well defined, for almost all paths the image is everywhere finite function on
edges, and the projection $\Phi$ is the canonical isomorphism of this space and
the boundary:
$$
(\Gamma(Sol_d,\mu), \hat \mu) = (F^d,\hat\mu). 
$$ 
here $\hat\mu \equiv \Phi(\mu^\infty)$
} 
\end{theorem}

\bigskip

Thus the Poisson--Furstenberg boundary of the pair $\Gamma (Sol_d, \mu)$
is canonically isomorphic as a measure space to the space 
$(F^d, \hat\mu)$ 
of integer-valued functions on edges provided with the image measure.

\begin{proof}{
We use a general method which is in more detals  published in paper by
author \cite{V}  in which we had consider the boundary of the groups for the
case when so called stable normal form exits. However, we shall make use of the
stable equivalence instead of stable normal forms. Again it suffices to prove
that the limits of functions depending on stable coordinates exhaust the whole
space of $\sigma$-measurable (boundary) functions. Assuming that such function
exists and is independent on all stable functions, approximate it by cylindric,
i.e.~$n$-stable functions depending on stable coordinates of length less than
$n$ (see the theorem of section~3). By the same reasons as in general method,
such functions must depend only on the coordinates of ($\beta(\gamma)$) that
can change when continuing the path $\gamma$, and thus, since the walk is
non-recurrent, the distance of the corresponding edges $\beta$ from zero must
be greater than some constant depending on $n$ for the paths $\gamma$ from a
set of measure close enough to one also depending on the choice of $n$.
Choosing, as above, a sequence $n_k$ so that these sets  of the numbers of
coordinates on which depend the  approximating functions, do not intersect, we
obtain  a contradiction (if the functions are not constant) with asymptotic
independence of values of the functionals $\beta(\gamma)$ for edges $\beta$
which are far enough from each other. Thus there is no functions except
constants that are orthogonal to all stable functions. Hence the boundary
$\Gamma (Sol_d, \mu)$ is mapped to $(F(E^d;\bf Z), \Phi(\mu^\infty))$ since
the image of a path under this mapping depends only on stable functionals of
the path, and this mapping is an isomorphism.}  \end{proof}

\begin{remark}
{
The properties of the measure $\Phi(\mu^\infty)$, i.e.~of the canonical measure
on the boundary of the free solvable group, are of great interest. This is a measure
on the space of integer-valued functions on the edges of the lattice. Its
one-dimensional distributions (the values on one edge) can be expressed by the
differences of the Green function of the simple $d$-dimensional
($d \geq 3$) random walk at the end and at the beginning of the edge. However, the
author knows nothing on the correlations of this natural measure. This measure
seems to be more natural than a similar measure on configurations, i.e.~on the space
of integer-valued functions at the nodes of the lattice which arises in the
wreath product model (see below).
}
  \end{remark}

Let us apply this method to other groups. Consider the wreath product
$\Z^d \wr H\equiv G_d(H)$, where $H$ is a cyclic group of finite or infinite order
(now this wreath product is called the ``lamplighter group'', since the walk on 
this group is
the walk on the lattice with simultaneous random switching lamps in each node of the 
lattice). This is also a solvable group of level $2$ with $d+1$ generators. 
It is convenient
to represent it as the skew product of $\Z^d$ and the space $F_0(\Z^d;H)$ of all finite
configurations, i.e.~$H$-valued functions on the lattice. 
Therefore, this group is naturally represented as a factor group of the group
$Sol_{d+1}$. If $d \geq 3$, then {\it the boundary
  $\Gamma(G_d(H), \mu)$, where $\mu$ is the uniform measure on generators, is not trivial},
this was first noted in \cite{KV}.  Easy to see that the boundary can be mapped
onto the space  $F(\Z^d;H)$ of all configurations: this mapping is just the ``final''
configuration, i.e.~with each node 
it associates the element of the group $H$ that was ``switched on'' when the walk
visited this node for the last time. However, till now it has not been clear whether
this mapping is an isomorphism, i.e.~whether each function on the boundary depends only
on the final configuration. This question was reduced to some problem on ordinary
simple walks on the lattice, but its solution was not obtained. In case when $H$ is 
a semigroup, the answer is positive (\cite{Peres}), but the method does not apply to
a group.

 \begin{corollary}{
The boundary  $\Gamma(G_d(H), \mu)$ is isomorphic to the space
   $(F(\Z^d;H), \hat \mu)$, where $\hat \mu$ is the "final" measure on the
space of configurations.
} 
\end{corollary}

\begin{proof}{
One may use the above theorem for the free solvable group
$Sol_{d+1}$ and the fact that the wreath product is its factor group. But it is more
instructive to carry out a similar argument for the wreath product itself making use of the
described method. For simplicity, consider the case
$d=3$ and $H=\Z$,  and denote the wreath product
$\Z^3 \wr \Z$ by $G_3$. In this case the wreath product has four generators, the first
three of them commute, and the commutators of any element from the subgroup
$Z^3$ formed by these three generators 
commute with the fourth generator too. Let us identify words
in generators with paths on the lattice $E^4$. It follows from above that two paths
are $G_3$-equivalent if and only if the projections of their ends on the 
sublattice formed by the first three axes coincide, and the algebraic sum of multiplicities
of occurrences of each edge parallel to the fourth axis is the same for both paths.
Further argument exactly reproduces the proof of the theorem  for the free solvable group.
}
\end{proof}

The previous attempts to prove this theorem ran across the following difficulty:
reduction of the problem on the boundary of the group $G_d$ to the walk
on the group $\Z_d$ (and not $\Z_{d+1}$) obscures the specific role of the fourth
generator. This 
leads to need to investigate the conditional process and to prove non-triviality
of its tail sigma-algebra, and this requires estimations of complicated functionals of
trajectories. In the above argument the distinction between generators is seen very well.
In particular, {\it the projection of multiplicities of edges parallel
to the fourth axis on the sublattice 
formed by the first three generators is exactly the final
configuration} which was discussed above.

We leave aside other examples, but only note that our method reduces all problems
on boundaries of finitely generated groups to problems on sigma-algebras (namely, the
sigma-algebras of $G$-equivalent paths) for the classical random walks on lattices of 
rank equal to the number of generators of the group. But it is not always easy to describe
the $G$-equivalence of paths (e.g.~for the braid groups). This method also shows that to 
find a normal form is to choose one path from the equivalence class, and the problem
of a 
minimal normal form is {\it an isoperimetric-like problem: to find a contour of minimal
length in a given class of (closed) paths}. For the continuous Heisenberg group this
is the classical isoperimetric problem. In the general case these problems
were called isoholonomic
(see \cite{VG}) ---  these are problems on the minimal length of a curve, given
fixed values of some family of $1$-forms.

For the free solvable group a normal form of elements also can be described, but it is
much more productive to take as coordinates the above considered generators rather than
the original ones.
It turns out that the description of the boundary does not involve
infinite words and even cannot be interpreted in ``cylindric'' (with respect to 
$S^\N$) terms. We shall consider these  questions in details elsewhere.

In conclusion, we note that it is very important  to make exact calculation of
the main constants --- logarithmic volume -
$$v=\lim_{N \mapsto \infty} \frac{\log W_{\leq N}}{N},$$ 
escape - 
$$c= \lim_{N \mapsto \infty}\frac{E_{\mu^{*N}}L(g)}{N},$$ 
and entropy
-$$h=\lim_{N \mapsto \infty}\frac{H(\mu^{*N})}{N},$$ 
(where $W_{\leq N}$ is a set of the elements of the groups of the
minimal length less or equal to $N$,  $\mu^{*N}$ is the N-th convolution of the
uniforn measure $\mu$ on, $L(g)$ is the minmal length of the element $g$ in 
those generators, and $E_\nu$ -- is an expectation with respect to measure
$\nu$). Perhaps this calculation for the free solvable groups of the level two
--  seems to be rather difficult, as well as an explicit calculation of the
values of harmonic functions. It is not yet known whether  fundamental
inequality 

$$h \leq l \cdot v$$ 
(see \cite{V}) turns into equality or not.   

The techniques we have used to find the boundaries 
by means of the spaces of paths, and the very realization
of groups, have a wide range of applications: the same tools may be used for
calculation of other boundaries (f.e.Martin boundary),
enumeration of measures with given cocycle, central measures, etc.

\bigskip
\bigskip

 ACKNOLEDGEMENT.The author expresses his gratitude to the International
Schroedinger  Institute where this 
paper was finished  for possibilities to work  in the institute in the
framework of semester on representation theory, and to
Dr.Natalia Tsilevich for English translation of the first version of the paper.


\end{document}